\documentclass[12pt,a4paper, reqno]{amsart}
\usepackage{amsfonts,amsmath,amsthm,amssymb, amscd}
\usepackage[top=2cm, bottom=5cm, left=3cm, right=1cm]{geometry}

\newtheorem{theorem}{Theorem}
\newtheorem{lemma}{Lemma}

\newtheorem{conjecture}{Conjecture}

\title[Commutative subalgebras of Grassmann algebra]{A remark on commutative subalgebras of \\Grassmann algebra}

\author{Ho-Hon Leung}
\address{Department of Mathematical Sciences, UAEU, Al-Ain, United Arab Emirates}
\email{hohon.leung@uaeu.ac.ae}

\keywords {Grassmann algebra, exterior algebra, maximal commutative subalgebra}
\thanks{
Supported by the UAEU grants: StartUp Grant 2016: G00002235}

\begin{document}
\maketitle

\begin{abstract}
Let $n=4k+1$ and $k\geq 4$. We show that there exists maximal commutative subalgebras (with respect to inclusion) of dimension less than $3\cdot 2^{n-2}$.
\end{abstract}

\section{Introduction}  \label{section1}
The Grassmann algebra  (exterior algebra)  $G(n)$ over a field $F$ of characteristic different from two  is the following  finite dimensional associative algebra of rank $n$:
\begin{align}\label{gwbRG}
G(n)&=F[x_1,\ldots, x_n]/\langle{x_ix_j+x_jx_i\mid 1\leq i,j\leq n}\rangle_F.
\end{align}It is clear that $G(n)$ has dimension $2^n$. In this short note, the term {\it maximal commutative subalgebra} refers to commutative subalgebras of $G(n)$ which are maximal with respect to inclusion.

If $n$ is even, Domokos and Zubor \cite{Domokos} showed that every maximal commutative subalgebra has dimension $3\cdot 2^{n-2}$. But, in general, the structures of maximal commutative algebras of $G(n)$ are less clear when $n$ is odd. one of the conjectures posted by Domokos and Zubor \cite{Domokos} is the following: 

\begin{conjecture} (Domokos and Zubor \cite[Conjecture 7.3]{Domokos}) \label{conjecture}
If $n=4k+1$ and $A$ is a maximal commutative subalgebra of the Grassmann algebra $G(n)$, then $\dim (A)\geq 3\cdot 2^{n-2}$.
\end{conjecture}

In a recent paper written by V. Bovdi and the author \cite{Bovdi}, we showed that the construction of maximal commutative subalgebras of $G(n)$ is closely related to the construction of certain intersecting families of subsets of odd size in the power set of $n$. Furthermore, based on the classical Fano plane, we constructed examples of maximal commutative subalgebras of $G(n)$ for odd $n$. In particular, we showed the following result numerically, which is a negative answer to Conjecture \ref{conjecture} for some values of $n$:

\begin{lemma} (Bovdi and Leung \cite[Corollary 4]{Bovdi}) \label{lemma}
Let $n=4k+1$ and $17\leq k< 1000$. There exist maximal commutative subalgebras in the Grassmann algebra $G(n)$ with dimension less than $3\cdot 2^{n-2}$.
\end{lemma}

\noindent We conjectured that Lemma \ref{lemma} still holds for $n=4k+1$ and $k\geq 4$ (see \cite[Conjecture 5]{Bovdi}). In this note, we answer this conjecture affirmatively:

\begin{theorem} \label{theorem}
Let $n=4k+1$ and $k\geq 4$. There exist maximal commutative subalgebras in the Grassmann algebra $G(n)$ with dimension less than $3\cdot 2^{n-2}$.
\end{theorem}

\section{Proof of Theorem \ref{theorem}} \label{section}

\noindent Let $k$ be any positive integer such that $k\geq 2$. We look at the following quantities:
\begin{align}
\label{equation0.1}  C_1 :=& 7\cdot \binom{4k+2}{2k} +\binom{4k+2}{2k+3}, \\
\nonumber  C_2 :=& \binom{4k+2}{2k+5} +\binom{7}{1}\cdot \binom{4k+2}{2k+4} +\binom{7}{2}\cdot\binom{4k+2}{2k+3}+7\cdot\binom{4k+2}{2k+2}+28\cdot \binom{4k+2}{2k+1} \\\label{equation0.2} &+\binom{7}{5}\cdot \binom{4k+2}{2k},  \\
\label{equation0.3}  C_i:=&  \binom{4k+9}{i}\mbox{ for }i\geq 2k+7\mbox{ and }i\mbox{ is odd},\\
\label{equation0.31} C_3 :=& \sum_{\{4k+9\geq i\geq 2k+7, i\mbox{ is odd}\}}C_i.
\end{align}Let $Q_k$ be the following quantity:
\begin{align}
\label{equation0.5}  Q_k &:= \frac{C_1 +C_2 +C_3}{2^{4k+7}}.
\end{align}Based on the construction of maximal commutative subalgebras of $G(n)$ in the paper written by Bovdi and the author \cite[Section 5]{Bovdi}, Theorem \ref{theorem} is equivalent to the following theorem:

\begin{theorem} \label{theorem1}
Let $k\geq 2$. Then \[ Q_k<1.\]
\end{theorem}

\noindent We define the variable $A$ as follows:
\begin{align}
\label{equation3}  A&:= \frac{128\cdot 16^k \cdot \Big( \sqrt{\pi}\cdot\Gamma(2k+6) -6k\cdot \Gamma\Big( 2k+\frac{9}{2}\Big) -13\cdot \Gamma\Big( 2k+\frac{9}{2}\Big) \Big)}{\sqrt{\pi} \cdot(2k+5)!}
\end{align}where $\Gamma(z)$ is the Gamma function. By the computer program Maple, it is shown that 
\begin{align}
\label{equation0.4}  C_3 = A. 
\end{align}

\noindent We note that the well-known Gamma function $\Gamma(z)$ has the following property: \[\Gamma(z+1)=z\Gamma(z)\]where $z$ is any complex number in the complex plane. Let $k$ be any positive integer. We expand $\Gamma(2k+4.5)$ as follows:
\begin{align}
\nonumber \Gamma\Big(2k+\frac{9}{2}\Big) &= \Big( 2k+\frac{7}{2}\Big) \Big( 2k+\frac{5}{2}\Big) \cdots \Big(\frac{1}{2}\Big)\Gamma\Big( \frac{1}{2}\Big) \\
\nonumber &= \Big( 2k+\frac{7}{2}\Big) \Big( 2k+\frac{5}{2}\Big) \cdots \Big(\frac{1}{2}\Big) \sqrt{\pi} \\
\nonumber &= \frac{1}{2^{2k+4}}\cdot (4k+7)(4k+5)\cdots 5\cdot 3\cdot 1\cdot \sqrt{\pi} \\
\nonumber &= \frac{(4k+7)!}{2^{2k+4}\cdot (4k+6)(4k+4)\cdots 6\cdot 4\cdot 2} \cdot \sqrt{\pi}\\
\label{equation1} &= \frac{(4k+7)!}{2^{4k+7} \cdot(2k+3)!} \cdot \sqrt{\pi}.
\end{align}By (\ref{equation1}), 
\begin{align}
\nonumber &128\cdot 16^k \cdot \Big( \sqrt{\pi}\cdot\Gamma(2k+6) -6k\cdot \Gamma\Big( 2k+\frac{9}{2}\Big) -13\cdot \Gamma\Big( 2k+\frac{9}{2}\Big) \Big) \\
\nonumber &= 2^{4k+7} \cdot \Big( \sqrt{\pi}\cdot (2k+5)! -6k\cdot \frac{(4k+7)!}{2^{4k+7}\cdot (2k+3)!}\cdot \sqrt{\pi} -13\cdot \frac{(4k+7)!}{2^{4k+7}\cdot (2k+3)!}\cdot \sqrt{\pi}  \Big)  \\
\label{equation2}  &= \sqrt{\pi}\cdot\Big( 2^{4k+7}\cdot (2k+5)! -6k\cdot\frac{(4k+7)!}{(2k+3)!}-13\cdot\frac{(4k+7)!}{(2k+3)!}  \Big).
\end{align}By (\ref{equation2}), we simplify the expression of $A$ in equation (\ref{equation3}) as follows:
\begin{align}
\nonumber  A&= 2^{4k+7} -\frac{(6k+13)\cdot (4k+7)!}{(2k+5)! \cdot (2k+3)!}\\
\label{equation4}  &= 2^{4k+7} -\frac{6k+13}{2k+5} \cdot \binom{4k+7}{2k+3}.
\end{align}By (\ref{equation0.1}), (\ref{equation0.2}), (\ref{equation0.5}), (\ref{equation3}), (\ref{equation0.4}), (\ref{equation4}), we write the expression $Q_k$ as follows:
\begin{align}
\label{equation5}  Q_k &=1+ \frac{D}{2^{4k+7}} -\frac{E}{2^{4k+7}}
\end{align}where the variables $D$ and $E$ are defined by 
\begin{align}
\label{equation6}  D &:= C_1+C_2, \\
\label{equation7}  E &:= \frac{6k+13}{2k+5}\cdot\binom{4k+7}{2k+3}.
\end{align}Theorem \ref{theorem1} is equivalent to the following theorem:
\begin{theorem} \label{theorem2}
Let $k$ be any positive integer. Then \[D<E.\]
\end{theorem}

\noindent We simplify equation (\ref{equation6}) as follows:
\begin{align}
\label{equation8}  D &= 35\cdot\binom{4k+2}{2k}+22\cdot\binom{4k+2}{2k+3} +\binom{4k+2}{2k+5}+7\cdot\binom{4k+2}{2k+4}+28\cdot \binom{4k+2}{2k+1}.
\end{align}We do the following algebraic manipulations on $D$ which will be needed later:
\begin{align}
\nonumber D\cdot\frac{(2k)!}{(4k+2)!}=&\frac{35}{(2k+2)!}+\frac{22\cdot 2k}{(2k+3)!}+\frac{(2k)(2k-1)(2k-2)}{(2k+5)!}+\frac{7(2k)(2k-1)}{(2k+4)!}\\\label{equation9}   &+\frac{28}{(2k+1)! \cdot (2k+1)}.
\end{align}
\begin{align}
\label{equation10} D\cdot\frac{(2k)!}{(4k+2)!}\cdot (2k+5)! &= \begin{aligned}[t] 
   35(2k+5)(2k+4)(2k+3)+22(2k)(2k+4)(2k+3)\\
+(2k)(2k-1)(2k-2)+7(2k)(2k-1)(2k+5)\\+\frac{28(2k+5)(2k+4)(2k+3)(2k+2)}{(2k+1)}.
        \end{aligned}
\end{align}Similarly, we have the following expression for $E$, 
\begin{align}
\label{equation11} E\cdot\frac{(2k)!}{(4k+2)!}\cdot (2k+5)! &= \frac{(6k+13)(4k+7)(4k+6)(4k+5)(4k+4)(4k+3)}{(2k+3)(2k+2)(2k+1)}.
\end{align}We multiply $(2k+1)(2k+2)(2k+3)$ to (\ref{equation10}) and (\ref{equation11}). The R.H.S. of these two equations become the following two expressions respectively: 
\begin{align}
\label{equation12}   \begin{aligned}[t] 
   35(2k+5)(2k+4)(2k+3)(2k+1)(2k+2)(2k+3)+22(2k)(2k+4)(2k+3)(2k+1)(2k+2)(2k+3)\\+(2k)(2k-1)(2k-2)(2k+1)(2k+2)(2k+3)
+7(2k)(2k-1)(2k+5)(2k+1)(2k+2)(2k+3)\\+28(2k+5)(2k+4)(2k+3)(2k+2)(2k+3)(2k+2),
        \end{aligned}
\end{align}and
\begin{align}
\label{equation13}   (6k+13)(4k+7)(4k+6)(4k+5)(4k+4)(4k+3).
\end{align}Let $D'$ and $E'$ be the expressions in (\ref{equation12}) and (\ref{equation13}) respectively. As polynomials in $k$, the dominating terms of $D'$ and $E'$ are $93\cdot 2^6\cdot k^6$ and $96\cdot 2^6\cdot k^6$ respectively. Hence, it is clear that \[D'< E'\mbox{ as }k\rightarrow \infty.\]More precisely, we obtain a surprising factorization for $E'-D'$ (by the help of Maple):
\begin{align}
\label{equation14}  E' -D' &= 24k(2k+5)(k+2)(k+1)(2k+3)^2.
\end{align}It is clear that $D' < E'$ for all positive integers $k$.

\noindent Finally, we note that the inequality $D' <E'$ is equivalent to the inequality $D<E$. Hence, Theorem \ref{theorem2} (which is equivalent to Theorem \ref{theorem1} and hence also Theorem \ref{theorem}) is proved. 

\section{Acknowledgement}
The idea of proving Theorem \ref{theorem1} by the use of the computer program Maple is due to Thotsaporn Thanatipanonda. The author is grateful that he shared his idea generously when the author visited him in Mahidol University International College in Bangkok in Summer 2018. Certainly, he takes the credit of providing the proof of Theorem \ref{theorem1}.

\end{document}